\documentclass[12pt]{amsart}
\usepackage{latexsym,fancyhdr,amssymb,color,amsmath,amsthm,graphicx,listings,comment}
\usepackage[section]{placeins}
\let\paragraph\subsection

\title{Interacting Geodesics on Discrete Manifolds}
\author{Oliver Knill}
\date{Mai 29, 2025}

\address{Department of Mathematics \\ Harvard University \\ Cambridge, MA, 02138 }


\begin{document}
\maketitle

\begin{abstract}
We define an evolution of multiple particles on a discrete manifold $G$.
Each particle alone moves on geodesics and particles can interact if they are
on the same facet. They move deterministically and reversibly on the 
frame bundle $P$ of the abstract simplicial complex $G$. 
Particles are signed and each is represented by a 
totally ordered maximal simplex $p \in P$ in $G$. 
The motion of divisors on $P$ also defines a time dependent reversible 
deformation of space. 
\end{abstract}

\section{The model}

\paragraph{}
Let $G$ be a {\bf finite abstract simplicial complex} of dimension $q$. 
This means that $G$ is a finite set of non-empty sets closed under 
the operation of taking finite non-empty subsets. We assume that the maximal dimension
$q={\rm max}_{x \in G} |x|-1={\rm max}_{x \in G} {\rm dim}(x)$ is positive. 
We also assume that $G$ is generated by its $q$-simplices
and that it has the property that the {\bf star} $U(z)=\{ x \in G, z \subset x\}$ of every 
$(q-1)$-simplex $z$ contains either one or two $q$-simplices. The $(q-1)$-simplices with $|U(x)|=2$ 
are called {\bf boundary walls}, the ones with $|U(x)|=3$ are {\bf interior walls}. Such simplicial
complexes generalize $q$-manifolds. For the later, every {\bf unit sphere} $S(x)=\overline{U(x)} \setminus U(x)$ of $G$
is required to be a $(q-1)$-sphere or a $(q-1)$-ball. The {\bf Euler characteristic} of $G$ is defined as 
$\chi(G) = \sum_{x \in G} \omega(x)$, where $\omega(x) = (-1)^{{\rm dim}(x)}$. 

\paragraph{}
The concept of manifold can be generalized a bit. We give now an example but only look at the 
boundary-less case. Inductively, a complex can be called a {\bf Dehn Sommerville sphere} 
of dimension $q$ if $\chi(G)=1+(-1)^q$ and every unit sphere $S(x)$ of $G$ is a Dehn-Sommerville sphere of
dimension $q-1$. The foundation of this inductive definition is that the 
empty complex $0$ is the $(-1)$-dimensional Dehn-Sommerville sphere. 
A simplicial complex $G$ is a {\bf q-Dehn-Sommerville manifold}, if every unit sphere
$S(x)$ is a $(q-1)$-Dehn-Sommerville sphere. For more details, see \cite{dehnsommervillegaussbonnet}.
Because the intersection of $q$ unit spheres is a $0$-Dehn-Sommerville sphere which agrees with a 0-sphere,
Dehn-Sommerville manifolds satisfy the assumptions we have asked for to define a geodesic flow.

\paragraph{}
Spheres form a submonoid of the monoid of simplicial complexes with {\bf join operation}
$A \oplus B = A \cup B \cup \{ a \cup b, a \in A, b \in B \}$. 
This can be seen by induction on dimension using the unit sphere formulas
$S_{A \oplus B}(x) = S_A(x) \oplus B$ for $x \in A$ and 
$S_{A \oplus B}(x) = A \oplus S_B(x)$ for $x \in B$ as well as the fact that
the {\bf simplex generating function} $f_G(t) = 1+\sum_{x \in G} t^{|x|}$ multiplies $f_{A \oplus B}(t)=f_A(t) f_B(t)$
and that $\chi(G)=1-f_G(-1)$, implying the {\bf reduced Euler characteristic} $1-\chi(G)$ is multiplicative.
Dehn-Sommerville spheres form a larger monoid, containing the sphere monoid, but unlike spheres,
Dehn-Sommerville spheres are recognizable by their combinatorial data $f_G(t)$ alone. As in a manifold,
also in a Dehn-Sommerville manifold, every unit sphere $S(x)$ is a join of two Dehn-Sommerville spheres 
$S^-(x)=\overline{ \{x\}} \setminus \{x\}=\{y, y \subset x, y \neq x\}$ 
and $S^+(x)=\overline{U(x)} \setminus \overline{\{x\}} = \{y, x \subset y, y \neq x\}$,
producing a hyperbolic structure $S(x) = S^-(x) \oplus S^+(x)$. 

\paragraph{}
The set $F$ of maximal simplices in $G$, the q-simplices, is also known as the set of {\bf facets}. The assumption
that $G$ is ``pure" is equivalent to the statement $G=\overline{F}$, where $\overline{F}$ is 
the {\bf closure} of $F$ in the finite non-Hausdorff {\bf Alexandrov topology} on $G$ \cite{Alexandroff1937,FiniteTopology},
where the open sets are the sets which are unions of stars $U(x)$ and the closed sets are 
sub-simplicial complexes of $G$. The stars form the vertices in the {\bf dual graph} $\hat{G}$, where 
two facets are connected, if their intersection is a $(q-1)$ simplex. For a $q$-manifold, 
$\hat{G}$ is a vertex-regular, triangle-free graph of degree $q+1$. 
It determines $G$ and could be completed to become a Dehn-Sommerville cell complex
dual to $G$. We call an {\bf ordered} maximal simplex $x=(x_0, \dots, x_q) \in F$ 
a {\bf frame}, because if the first element $x_0$ is considered to be the {\bf base point} then the $q$ 
edges $(x_0,x_k), k=1, \dots, q$ form the analogue of a {\bf coordinate system}.
The set $P$ of all frames forms a {\bf discrete frame bundle} over $F$. 
It is a {\bf principle $S_{q+1}$-fiber bundle} $\pi: P \to F$ with the symmetric group $S_{q+1}$ as 
structure group and $|P|=f_q (q+1)!$ elements. An individual {\bf particle} is a point in $P$ with an additional sign.
The {\bf particle position} of $x \in P$ is $\pi(x) \in F$. The orientation in $S_{q+1}$ can be thought of 
as representing a {\bf particle momentum}. We write $-(x_0,\dots, x_q)$ for a negative particle.

\paragraph{}
Every ordered $q$-simplex $x=(x_0, \dots, x_q)=(x_0,z) \in G$ defines a {\bf base vertex} $x_0$
as well as a {\bf wall} $z$, a $(q-1)$ simplex. The assumption for $G$ gives now two possibilities: 
either $U(z)=\{ (x_0,z), (x_0,z'),z \}$ or $U(x)=\{ (x_0,z),z\}$, so that $x \in P$ has a partner 
$x'=(x_0', \dots, x_q)$. There are two possibilities. 
Either $x_0'$ is different to $x_0$ which means that 
$U(z) = U( (x_1, \dots, x_q) ) = \{ z,x,x' \}$ or then $x_0'=x_0$ so that
$U(z) = \{ z,x \}$. On a $q$-manifold without boundary, only the first case happens. 
Also on a Dehn-Sommerville manifold, only the first case happens, the reason being that
the only Dehn-Sommerville sphere of dimension $0$ is the usual $0$-dimensional sphere
consisting of two disjoint points. One could look at Dehn-Sommerville manifolds with boundary,
but we do not do that here. 

\paragraph{}
{\bf Particles} are oriented $q$-simplices that are additionally equipped with a 
{\bf signature} $+$ or $-$. This allows to describe a {\bf particle configuration}, a finite collection
of particles, as a {\bf divisor} on $P$. The language of divisors comes from discrete Riemann-Roch theory 
\cite{BakerNorine2007}. There are other reasons why we like to consider a signature for particles:
the most important one is that it allows to write the geodesic motion and also the 
$n$-particle evolution $T$ as $T=BA$, where $B,A$ are both {\bf involutions}. This has 
a practical advantage, because the {\bf backwards evolution} $T^{-1}=AB$ has the same structure. 
By the way, also the translation group and rotation group of Euclidean space can be Coxeterized as such:
the composition of two point reflections is a translation, the composition of two hyperplane 
reflections is a rotation, which is the point of view of ``Spiegelungsgeometrie" 
\cite{Bachmann1959,JegerTransformationGeometry}.

\paragraph{}
Also physics motivates to consider signed particles: we observe signatures in nature in the form of {\bf mass}, 
{\bf charge}, {\bf spin} or {\bf polarization}. Modern quantum field theory is built
on creation and annihilation processes, or pictures seeing virtual particle-anti-particle pairs filling out space.
In mathematics, arithmetic is a strong motivator for positive and negative geometry parts 
\cite{numbersandgraphs}. One often wants to complete a monoid to a group, by introducing negative or
inverse elements. Simplicial complexes which form a monoid with respect to disjoint union can
be group-completed as such and then together with the Shannon product produce a ring, the {\bf Shannon ring}. 
Finally, we can think of adding negative particles as a completion fixing one of the 
asymmetries in the definition of $T$: we can either left or right rotate in the definition of 
geodesic flow. We will also see that negative particles can be associated with $q$ identical 
positive particles on the same $x \in P$. 

\paragraph{}
A {\bf single particle} is a positive or negative element in $P$. Its evolution will move on a
{\bf geodesic}. A finite particle configuration $(X_1, \dots, X_{n^+}, Y_1, \dots, Y_{n^-}$ 
is a collection of finitely many particles, where the $X_i$ enumerate the 
{\bf positive particles} and the $Y_i$ enumerate the {\bf negative particles}. 
Particles form a {\bf divisor} on $P$ that is an element in the {\bf free Abelian group} on $P$. 
Also the particle positions $\pi(X_i),\pi(Y_i) \in F$ define a {\bf divisor} $X$ on $F$, if 
$\pi:P \to F$ is the projection, which simply forgets about the total order structure. 
We could define ${\rm deg}_x(X) = \sum_{p \in \pi^{-1}(x)} X(p)$ 
counting how many particles there are on a facet $x \in F$. It is even possible to
see a particle configuration as a $Z_\{q+1\}$ bundle over $P$, as we will just see. 

\paragraph{}
On the bundle $P$, there can be several particles on 
the same point $x \in P$. Having $q$ positive
particles on the spot $x \in P$ is the same than having $1$ negative particle there. 
The reason is that $q+1$ positive or negative particles on the same $x \in P$ do not influence the
rest of the particles. They form a part of the particle cloud that are bounced around the 
other particles but they do not influence the rest. We call a group of $q+1$ particles on the same spot
an {\bf eddie}, similarly as classical mechanics, where the motion of a test particle in a 
time-dependent gravitational field of other bodies is considered. In classical mechanics such systems are
called {\bf restricted n-body problems} and motivate to look at time dependent Hamiltonians $H(p,q,t)$, 
defining via the Hamilton equation $q'(t)=H_p(q,p,t),p'=-H_q(q,p,t)$ a symplectic curve $\phi_t$ in the
cotangent bundle $T^*M$ or a manifold $M$. Such flows play an important role in symplectic geometry
\cite{HoferZehnder1994}. 

\paragraph{}
Since $q+1$ particles at the same point $x \in P$ form an
``invisible" part of space, we can ignore it when experimenting with this game. 
One benefit to ignore such ``ether" components and assume to have maximally $q^{|P|}$ possible 
particle configurations. But the number of possibilities is still large.
On the smallest 2-manifold, the octahedron $G=K_{2,2,2}$, 
where $|P|=6*3!=36$, we have already $3^{36}=1.5 \cdot 10^{17}$ possible particle configurations.
An advantage of seeing an eddie is to make the particle stream visible. It is a tracer that
does not influence the motion of the particles. Similarly as in Vlasov dynamics, where the eddie motion 
is a deformation of a manifold equipped with a probability measure $\mu$ 
given by $f''=\int_M F(f(x),f(y)) d\mu(y)$ which defining a symplectic map $(f(t),g(t))$ 
on the cotangent bundle $T^*M$ (see \cite{Vlasov,BrHe77,Vlasovgasboundary}).
We can see ``space" as part of ``matter" if one sees space built up by particle-antiparticle pairs on 
each point $p \in P$ or then by seeing space as a configuration, where every point in $P$ is 
occupied by one particle. In the former case, the dynamics $T$ will be an involution, in the later,
it is the geodesic flow for one particle. 

\paragraph{}
Every particle configuration defines
a time-dependent volume preserving transformation $T_t$ on $P$, 
meaning that the time evolution defines a path in the large
permutation group $S_{|P|}$ of $P$. It has the property that $T_t(x),T_{t+1}(x)$ are 
adjacent simplices for every $x \in P$ and all $t \in \mathbb{Z}$. 
Getting such a generalized dynamical system has been the original motivation to look 
at particle dynamics. The geodesic flow on $G$ defined in \cite{knillgeodesics} describing the 
motion of a single particle is a defined by local rules so that we can not expect that any two
pairs of points $x,y$ can be connected by a geodesic. We can hope however that there exists a 
small time interval $I=[0,\tau]$ with $\tau$ logarithmically small with respect to $|G|$ such 
that if we start an eddie at $x$ and pick a point $y \in F$, there is a time $t \in I$ such 
that $T^n x = T_{t+n} \cdots T_{t+1} T_{t} x = y$. We expect this to happen if space is filled 
with a sufficiently dense particle gas. 

\paragraph{}
Now, lets finally get to the definition of the {\bf evolution map} $T=BA$. 
The goal is to get a permutation on all possible particle configurations. We write $T$ as
a product of 2 involutions: $A(x)=-x'$ and $B(+x)=-L^{k-l}(x)$ and 
$B(-x)=+L^{l-k}(x)=+R^{k-l}(x)$, where $R(x_0, \dots, x_q)=(x_1, \dots, x_q,x_0)$ 
is {\bf right rotation} and $L(x_0, \dots, x_q)=(x_q,x_0, \dots, x_{q-1})$ is 
{\bf left rotation} and where $k,l$ are the number of positive respectively negative
particles on $x$ which project to the same {\bf position} $\pi(x) \in F$. 
For a single particle $n=1$, we get the geodesic flow as defined in 
\cite{knillgeodesics}. Note that because $T=BA$, we have $T^{-1}=AB$. 
We insist on particles to interact locally at the same point in $F$ only. This is
justified by {\bf causality}: there should not be any instantaneous force between parts 
of space which which have a positive distance. 

\paragraph{}
A positive particle $X_i$ and negative particle $Y_j$ at the same point $x \in F$
can not be ignored in general as they can have different orientation. We can for example
have $X_1=\{1,2,3\}$ and $Y_1=-\{1,3,2\}$ in a $2$-dimensional complex. The particle configuration
$\{ X_1,Y_1 \}$ produces a ``blinker", a periodic point of period $2$. If now a third particle
hits this blinker in the right time so that all three are together at the same $x \in F$, 
then we have split the molecule and got three particles moving away on geodesics. 
But if they occupy the same point in $P$, then we can ignore such a pair. When being hit
by a third particle, they move but then again remain a blinker. In any way, 
they do not contribute to the motion of other particles. But they can be 
bounced around like a pollen particle in Brownian motion. 

\paragraph{}
As pointed out already, having $q+1$ particles of the same sign at the same point
can be ignored if they have the same orientation. More generally, a
set of particles $X_i,Y_i$ with the same orientation and for which the
total degree is $0$ modulo $q+1$ can be removed as they do not affect the motion 
of the other particles.  The total number of 
possible relevant particle configurations we would have to consider 
therefore is $f_q (q+1)! (q+1)$. We
can also look just at $q$ identical positive particles as one negative 
particle or $q$ identical negative particles as one positive particle. 

\paragraph{}
There is an other extension which can be motivated both from mathematics as well 
as by physics, especially when thinking about the rotation in the fiber as spin and 
by trying to implement ``spin 1/2".
Instead of taking the rotation group $Z_{q+1}$ on the fiber as the
core of the interaction, we can generalize the fiber action to become
the {\bf dihedral group} $D_{q+1}$. It is more natural in the sense of
\cite{GraphsGroupsGeometry} as the group structure on this group can be derived from a 
metric on its set. An element in $D_{q+1}$ contains now both involutions and
and rotations. A rotation by $1$ be written as the product of two reflections, a 
reflection at the first element $1$ and a reflection at the last element. A Fermionic
evolution is given by $T=A C$, where $C$ is a reflection, if the number of
Fermions present at $\pi(x)$ is odd and the identity otherwise. We have a 
{\bf Pauli principle} in that already 2 Fermion particles on the same point 
$x \in P$ can be ignored: they do not influence the rest of the particles and 
remain together.

\paragraph{}
One goal has been to get an evolution in which we have no free parameters.
Large language models with billions of parameters
have demonstrated, that rather convincing intelligence can be faked 
using with cleverly trained perceptron and attention layers, which then boils down to 
some randomized data fitting. There is little real understanding however. For example, if a LLM 
draws pictures of turbulence it might not yet see the underlying partial differential equations or
have insight about conservation laws. According to Freeman Dyson, Enrico Fermi once told that Jonny Von Neumann 
has pointed out that with ``4 parameters one can fit an elephant" and that ``with 5, we can even make 
it wiggle its trunk".
Mathematicians like inevitability: examples are singling out {\bf Euler characteristic} as the
only valuation that is invariant under Barycentric refinements or to single out {\bf entropy} as 
the only functional on finite probability spaces satisfying some arithmetic properties among
probability spaces \cite{Shannon48}. In any way, we feel that real understanding of a phenomenon can
also be measures in how much assumptions and free parameters one has to assume. Real understanding would
mean we can compute numbers like the masses of particles or the fine structure constant from basic principles.

\paragraph{}
To summarize: \\

{\bf 1)} We have defined a {|bf global reversible dynamical system} that describes particles on rather general
classes of simplicial complexes $G$. Particles are divisors on the frame bundle $P$ of the facets $F$ in $G$.
The frame work can be used complexes that satisfy rather weak manifold structure, like 
Dehn-Sommerville manifolds. \\
{\bf 2)} {\bf Single particles move on geodesics}.
The total {\bf number of particles} $\sum_x X(x) + Y(x)$ is preserved. 
Also the total {\bf degree} $\sum_x X(x)-Y(x)$ is preserved.  \\
{\bf 3)} The particle motion is a {\bf reversible cellular automaton} \cite{Hed69}.
All signals propagate with a universal maximal ``speed of light". Changing the position or momentum of
a particle at some point $x \in F$ can not affect particles in distance $n$ 
in less than $n$ time steps.  \\
{\bf 4)} If $X(x)=0 \; {\rm mod} (q+1)$ for all $x \in P$, then $T^2(X)=X$.
An example is to take a particle and anti-particle at every $p \in P$.
If particles cover space as such their motion encodes the 
geometry $G$ by defining a graph on $P$. Two points $p,q$ are connected if $T(X)=Y$. 
An other way to see {\bf space} as part of the {\bf matter} frame work is to assume $X(x)=1$ for every $x \in P$. 
In this case, space itself deforms under {\bf time} evolution and the time evolution also reveals the 
structure of space. The game is a concrete allegory for an ancient story \cite{Weyl1950} (written in 1918). \\
{\bf 5)} The particle motion defines a time-dependent permutation $T_t$ on $P$ which describes the 
{\bf eddie motion} of particles that do not influence the others, like for example a particle-anti-particle pair. 
Unlike the geodesic dynamics $T$ for a single particle which is deterministic and defines a permutation of $P$, 
the time dependent permutation can be much richer. We expect with relatively low particle density to get a
time-dependent permutation $T_t$ which is transitive, meaning that every $x,y \in P$ can be connected in the eddie
motion like planets stirring space wround the sun can be used to shoot space proves rather arbitrarily by choosing
the timing cleverly. \\
{\bf 6)} We have a simple {\bf space-time-matter game} to play with. It was not designed to 
be a model for {\bf particle physics} but to have a toy model that {\bf maximizes interest} 
in the form of questions we can ask for it and {\bf minimizes the complexity} to build up. We feel that
this is a game that works. There are countless ambitious programs to conquer fundamental questions in physics
but some is mathematically not sound. 

\section{Experiments}

\paragraph{}
In general, if at every point $x \in P$, there exists exactly one positive
particle, then all particles just move on geodesics. The reason is that there the 
fiber above every $\pi(x) \in F$ has $(q+1)!$ points which is a multiple of $q+1$.
Space deforms just according to the geodesic flow equations. The same happens if at every point, 
we have a single negative particle is placed, leading to an essentially reversed flow. 
If at every point there is one positive and one negative particle, then the
configuration is a {\bf blinker}, a periodic point of period 2 in the lingo of Conway's game of life.
If at every point a constant number of $k$ particles are present, then the flow is
conjugated still to the standard geodesic flow as long as the greatest common divisor between 
$k$ and $q+1$ is ${\rm gcd}(k,q+1)=1$. 

\paragraph{}
If $q=1$, then $G$ is a finite union of disjoint path graphs or circular
graphs. Since particles stay on the same connected components, we can look
at each component separately. Particles live on edges. There are no positive
or negative particles in one dimension, just one type of particles.
Any pair of particles on the same point defines a blinker (again using the terminology of the 
game of life \cite{GardnerLife}). All other particles do move freely and
change direction at a pair collision where two particles hit common point $x \in F$. 

\paragraph{}
For example with 2 particles $X=(1,2)$  and $Y=(3,2)$ we have $T(X)=(2,3), T(Y)=(2,3)$. 
At that time, the second part $B$ of the dynamics $T=BA$ is just negation so that
$T(2,3)=(1,2)$ and $T(2,3)=(3,2)$. The two particles have bounced. If the timing is such 
that the two particles never occupy the same simplex, the particles just pass through each
other. For example with $X=(2,3)$ and $Y=(4,3)$, we have $T(X)=(3,4)$ and $T(Y)=(3,2)$.
The two particles have passed each other without interacting. 

\paragraph{}
For $q=2$, we can for the first time consider particles of different signature. 
We also observe moving {\bf molecules}, particles that move and stay together along 
geodesics. Being hit by a third particle can scatter them. We see some examples in the 
illustration. 

\begin{figure}[!htpb]
\scalebox{0.75}{\includegraphics{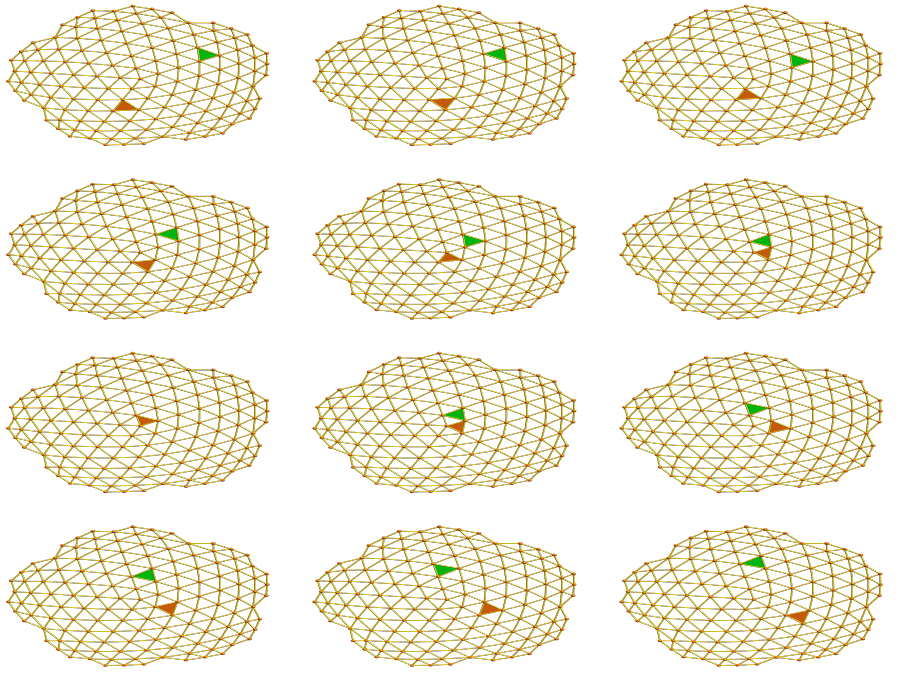}}
\label{scattering1}
\caption{
We see a positive and negative particle scatter in a 2 dimensional
complex. 
}
\end{figure}

\begin{figure}[!htpb]
\scalebox{0.75}{\includegraphics{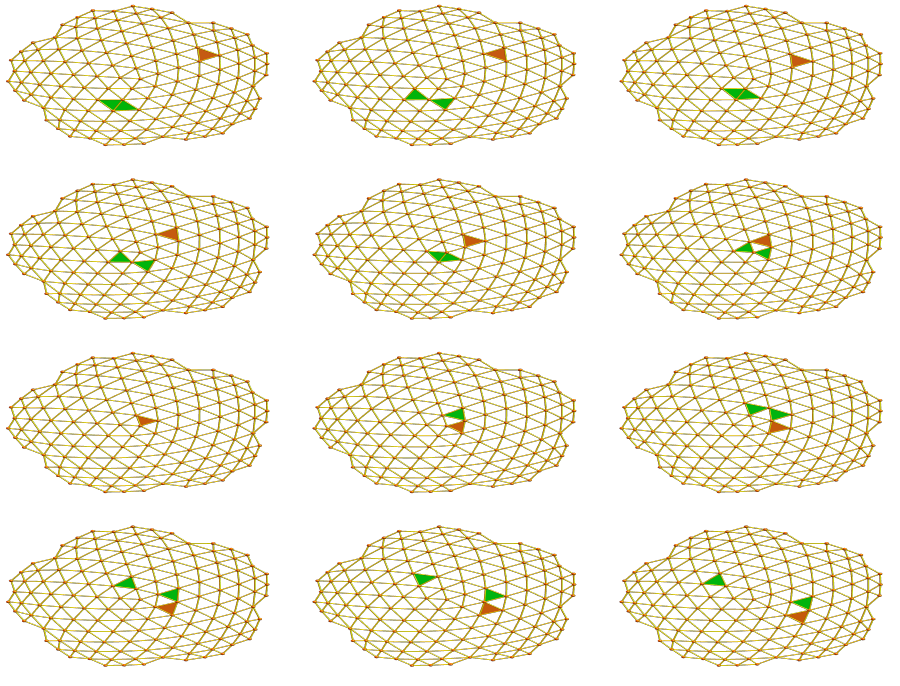}}
\label{scattering2}
\caption{
Here we see a molecule of two positive particles colliding with 
a negative particle. In this case, one of the positive parts gets
attached to the negative particle and moves away. 
}
\end{figure}

\begin{figure}[!htpb]
\scalebox{0.75}{\includegraphics{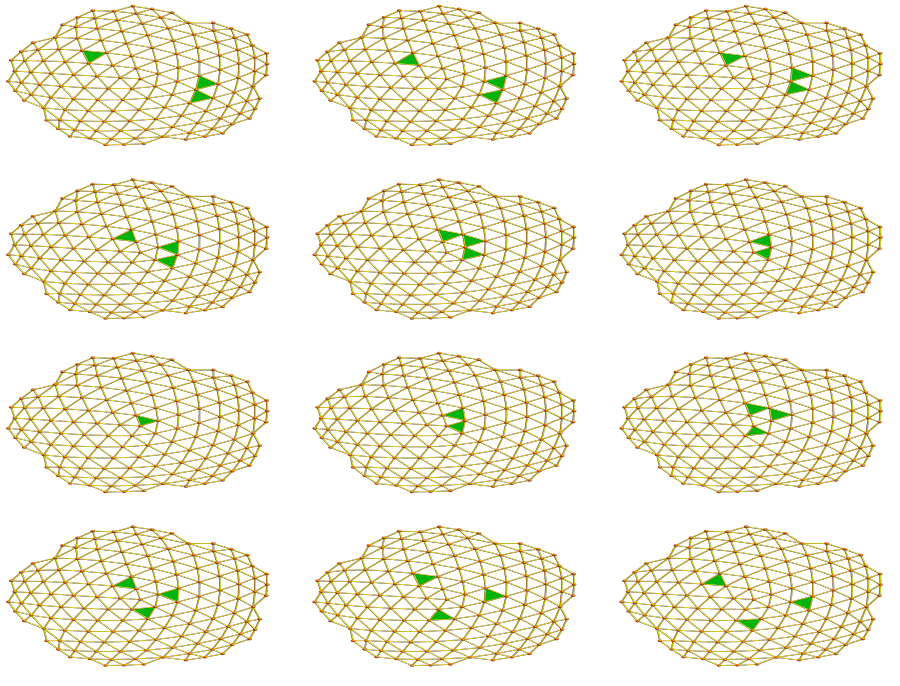}}
\label{scattering3}
\caption{
Three positive particles, where one is a molecule scatter and 
leave all individually. 
}
\end{figure}

\begin{figure}[!htpb]
\scalebox{0.75}{\includegraphics{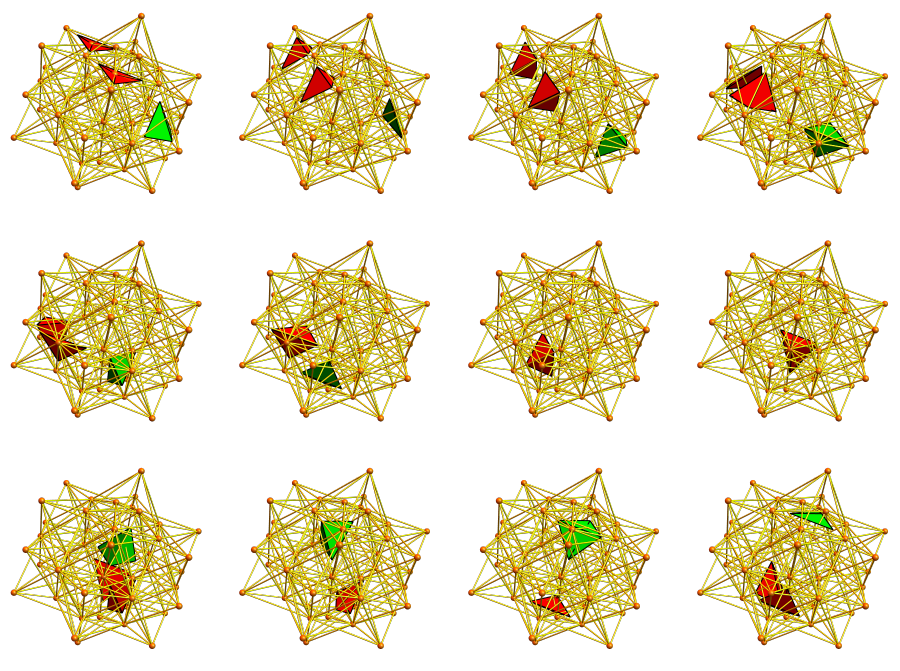}}
\label{scattering4}
\caption{
Three particles scattering in 3 dimensions. The 
method and code works without modification in any dimension. 
}
\end{figure}

\begin{figure}[!htpb]
\scalebox{0.75}{\includegraphics{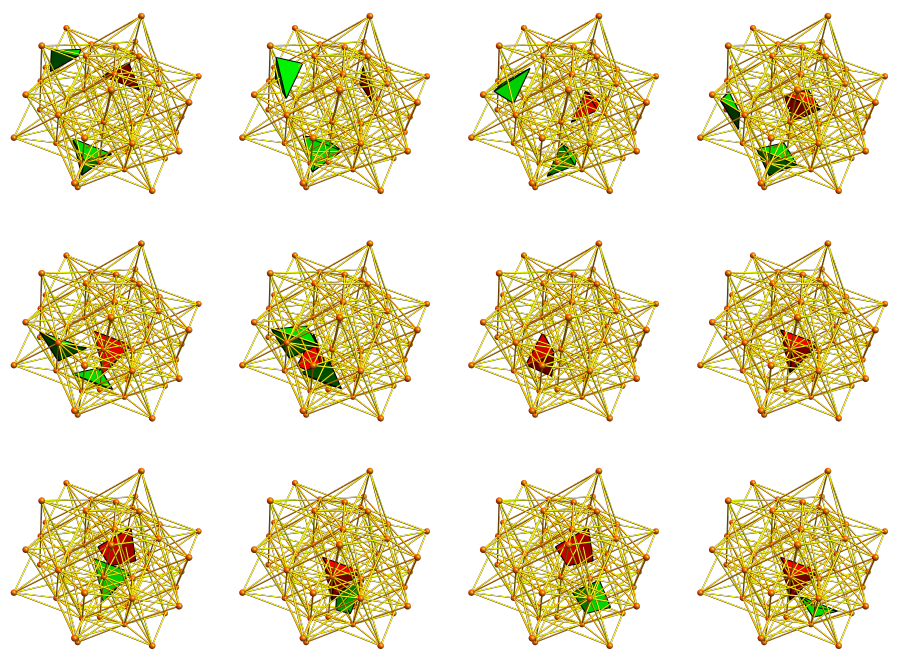}}
\label{scattering5}
\caption{
Two positive and one negative particles form a triple collision,
leaving a particle-anti-particle blinker and one freely moving particle. 
}
\end{figure}

\section{Remarks}

\paragraph{}
The definition given here was motivated by classical mechanics, noting that 
particles interact by transferring momentum. To have geodesics at the core of the definition
is motivated of course by physics and especially by relativity. 
By making the interaction local and space independent, the dynamics is a cellular 
automaton (CA). The proposition to use CA as the foundation of physics has been suggested
several times. It started with Ulam and von Neumann in 1940 and is part of ambitious programs like
\cite{Wolfram2002, tHooftCA}. One of the issues with these suggestions is
that they appear still too be too rich. Cellular automata can naturally be seen as discrete analogues of 
{\bf partial differential equations}, where both time, space as well as well as the function values 
are discrete or even finite and where the rules of the time propagation are the same at every cell. 
For more about cellular automata, see \cite{Wolfram86,SchiffCA,ChopardDrozCA,IlachinskiCA,HawkinsCA}.

\paragraph{}
We need an Abelian interaction rule if we do not want to have the interaction depend on the order in which 
particles arrive at the same point. We do not not discuss such an extension here.
The evolution $T$ should preserve the number of particles. 
Having an Abelian addition in the fiber appears to be the simplest interaction 
that taps directly into the definition of the 
geodesic flow. While particles move alone, they can be stuck temporarily and
revived later when an other particle moves through. Let $x$ for example be a boundary 
facet and assume it is the place where a positive and negative particle meet. The
two particles are now stuck there. But if an other single geodesic hits it, we have
three particles moving away, possibly in different directions. 

\paragraph{}
We can also look at space itself as a particle configuration and the dynamics given by
this particle configuration, to encode the topology. 
Look at a particle configuration in which every point $p \in P$ contains
a particle as well as an anti-particle. Each of these pairs just bounces forth
and back since $T(X)=A(X)$ with the involution $A$. This periodic dynamics encodes the
simplicial complex, because it tells which parts of space are glued together. The dual 
graph $\hat{G}$, telling which maximal simplices are attached to each other completely determines
the geometry of $G$ because the set of facets is $F$ is dense in $G$ in the Alexandrov topology: 
the way how $F$ are connected determines how the closure is connected: the closure of $F$ 
is just the set of all subsets of maximal simplices.

\paragraph{}
Seeing a particle configuration as generalized space should allow to see them
as elements in a delta set. (For abstract delta-sets in which we start with a finite set of sets
and a Dirac matrix $d+d^*$, see \cite{fusion2}.) We plan to work this out more in the future.
We want to see traditional space itself as the situation, where one particle is at every $x \in F$. 
This is motivated by discrete Riemann-Roch, where the canonical divisor plays the role of "space".
If $f_k(G)$ is the number of elements in $G$ of dimension $k$, we want that
$\chi(G) = \sum_{k=0}^q (-1)^k f_k(G)$ extends to $\chi(X) = \sum_{k=0}^q (-1)^k f_k(X)$
and that Gauss-Bonnet generalizes from spaces $G$ to particles $X$, divisors on $P$. 
The Levitt {\bf curvature of a vertex} $v \in G$
$ K(v) = 1-\sum_{k=0}^{d-1} \frac{(-1)^{k} f_k(S(v))}{k+1}$ \cite{Levitt1992,cherngaussbonnet}
should generalize from space $G$ to particles $X$, if $S(v)=\overline{U(v)} \setminus U(v)$.
The polynomial $f_X(t)=1+f_0(X) t + \cdots + f_d(X) t^{d+1}$ is the
{\bf simplex generating function} of $X$ satisfying $\chi(X)=1-f_X(-1)$.
The Gauss-Bonnet relation $\sum_v K(v) = \chi(G)$ should still be written as
$f_G'(t) = \sum_{v \in V} f_{S(v)}(t)$.

\paragraph{}
It is not clear whether the gas dynamics defined here can behave
like a discrete fluid if a large number of particles are taken. 
Even so it might not look so, Boltzmann models of particle dynamics
have been quite successful. There are also various cellular automata models.
The {\bf Hardy-Pomeau-Pazzis (HPP)} model is a lattice gas automaton on 
a two dimensional square lattice. It was proposed 50 years ago to 
simulate gases and liquids \cite{HaPo72,Har+76}. 
The Frisch-Hasslacher-Pomeau (FHP) model works on the hexagonal
lattice, but has a probabilistic interaction rule.
See also \cite{CardinalsChaos}.  Like any cellular automaton, 
it can also be simulated in almost periodic settings \cite{HofKnill}. 
Defining reasonable lattice gas cellular automata in higher dimensions would
in each lattice case require some choice of how collision rules work. 

\paragraph{}
Unlike particle systems which are based on potentials coming from Green functions,
we have defined here a globally defined dynamics which also works if two or more particles
interact. For particle interactions coming from a non-linear
smooth force, there is no danger coming from collisions but they are all non-relativistic
in the sense that particles in positive distance immediately interact. While no such issues
appear in the hard core interaction of balls, there is a difficulty for the hard core gas: one has to
neglect triple and higher collision because their outcome is undetermined.  
For Newton type interactions, where the potential comes from Green function,
already double collisions are not defined. Additionally, all models with pair interaction 
over some distance violates causality. For realistical relativistic models, particle 
forces exchange information via gauge Bosons. 

\paragraph{}
When looking for more general discrete particle models based on what we have looked at here,
one in principle has a lot of choices. 
One could define a map which maps $k$ points $p_i \in P$ with $\pi(p_i)=x \in F$
into an other $k$ tuple of points $T(p_i)$. In the case of $1$ point, the motion 
needs to produce the geodesic flow. We could make the interaction to depend on
the order in which the points are considered but that would require some choice.
Motivated by physics is to look at Fermion parts, where the rotation in the fiber is 
is generated by reflections. Then one has a Pauli principle and two identical Fermions 
at the same spot $p \in P$ do not matter. Such Fermions however stay localized and 
must be pushed around by Bosons. Their motion is much slower than the motion of the 
Bosons. Their speed will depend on the frequency with which they are visited by Bosons
doing the turns in and so the propagation in space. 


\section{Code}

\begin{tiny}
\lstset{language=Mathematica} \lstset{frameround=fttt}
\begin{lstlisting}[frame=single]
Generate[A_]:=If[A=={},{},Sort[Delete[Union[Sort[Flatten[Map[Subsets,Map[Sort,A]],1]]],1]]];
Whitney[s_]:=Map[Sort,Generate[FindClique[s,Infinity,All]]];L=Length;Co=Complement;Po=Position;
Facets[G_]:=Select[G,(L[#]==Max[Map[L,G]]  ) &];  RL=RotateLeft;  RR=RotateRight; S=Sort;
Bundle[G_]:=Module[{F=Facets[G]},Flatten[Table[Permutations[F[[k]]],{k,L[F]}],1]];  Se=Select;
OpenStar[G_,x_]:=Select[G,SubsetQ[#,x]&];         Stable[G_,x_]:=Complement[OpenStar[G,x],{x}];
M[z_]:=Module[{U=Co[Stable[G,z],{S[z]}]},Table[First[Co[U[[j]],z]],{j,L[U]}]];
a[x_]:=Module[{y=x,u},u=M[Delete[x,1]];If[L[u]==1,y[[1]]=u[[1]],y[[1]]=Co[u,{x[[1]]}][[1]]];y];
A[{X_,Y_}]:={Table[y=Y[[k]];a[y],{k,L[Y]}],Table[x=X[[k]];a[x],{k,L[X]}]};
B[{X_,Y_}]:={Table[y=Y[[k]];RL[y,L[Se[Map[S,X],S[y]==#&]]-L[Se[Map[S,Y],S[y]==#&]]],{k,L[Y]}],
             Table[x=X[[k]];RL[x,L[Se[Map[S,X],S[x]==#&]]-L[Se[Map[S,Y],S[x]==#&]]],{k,L[X]}]};
T[{X_,Y_}]:=B[A[{X,Y}]]; TI[{X_,Y_}]:=A[B[{X,Y}]]; Orbit[{X_,Y_},n_]:=NestList[T,{X,Y},n];
Orbit[{X_,Y_}]:=Module[{U={X,Y},o},o={U};While[U=T[U];Not[MemberQ[o,U]],o=Append[o,U]];o];
G=Whitney[CompleteGraph[{2,2,2,2}]];P=Bundle[G];{X,Y}={{P[[1]]},{P[[9]]}};o=Orbit[{X,Y}]; o
Q=Tuples[P,2];m=0;Do[u=L[Orbit[{{P[[1]],Q[[k,1]]},{Q[[k,2]]}}]];If[u>m,Print[u];m=u],{k,L[Q]}];


\end{lstlisting}
\end{tiny}

\section{Questions}

\paragraph{}
{\bf Q1:} The particle motion defines a time dependent permutation flow of $P$. We expect 
that this time dependent permutation is transitive over a small time interval already. 
Can we quantify this We expect that on an interval $I=[0, \log(|P|)]$ with $\log|P|<{\rm diam}(G)$,
the eddy dynamics becomes transitive: For all $x,y \in P$, there exists $t \in I$ such 
that $T^t*x=y$. We can even hope that on an eventually larger time interval like $\sqrt{|P|}$,
we can get $T^t(x)=y$ for $t \leq C d(x,y)$ for some constant $C(G) \to 1$ for $|G| \to \infty$.  \\

{\bf Q2:}  We would love to have a cohomology for particles which is compatible with Riemann-Roch.
This requires to write $X$ as a delta set, leading to face maps $d$ producing a matrix $D=d+d^*$
for which the kernels of the blocks of $D^2$ produce the Betti numbers. 
The Euler characteristic of a particle configuration $X$ should be 
$\chi(G)+{\rm deg}(X)$. We can look at the Laplacian of the dual graph 
$\hat{G}$ and have so a notion of linear equivalence. Is there a meaning for the
dimension function $l(X)$ which appears in Riemann Roch? How does this change?  \\

{\bf Q3:} What happens with particle configurations that are equivalent in the sense of 
Riemann Roch? If we take $q+1$ particles on a point $p \in P$ with different orientations
and evolve them, we get a {\bf chip firing situation}. A concrete question is:
Assume $X \sim Y$ are linearly equivalent divisors on $P$. Are then also 
$T(X),T(Y)$ linearly equivalent? This would mean that the Riemann-Roch quantity $l(X)$
entering the Riemann-Roch equation $l(X) - l(K-K)=\chi(G)+ {\rm deg}(X)$ is preserved 
by $T$. 

\bibliographystyle{plain}

\end{document}